\documentclass[12pt]{amsart}
\usepackage{amssymb, epsf, colordvi}
\usepackage{multirow}
\newif\ifEPSF
\EPSFtrue %
\numberwithin{equation}{section}

\newtheorem{thm}{Theorem}
\numberwithin{thm}{section}
\newtheorem{prop}[thm]{Proposition}
\newtheorem{lemma}[thm]{Lemma}
\newtheorem{cor}[thm]{Corollary}

\newtheorem{example}[thm]{Example}
\newtheorem{remark}[thm]{Remark}

\newenvironment{ex}{\begin{example}\rm}{\end{example}}
\newenvironment{rem}{\begin{remark}\rm}{\end{remark}}
\newcounter{FNC}[page]
\def\newfootnote#1{{\addtocounter{FNC}{2}$^\fnsymbol{FNC}$%
     \let\thefootnote\relax\footnotetext{$^\fnsymbol{FNC}$#1}}}

\newcommand{\EOP}{\qquad\raisebox{-1pt}{\epsfysize=8pt\epsffile{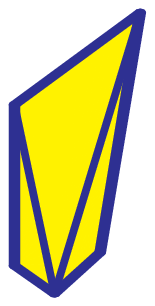}}\medskip}
\newcommand{\QED}{\ifEPSF \EOP \else \qquad$\Box$\medskip \fi}

\copyrightinfo{}{}
\newcommand{\calA}{\mathcal{A}}
\newcommand{\calB}{\mathcal{B}}
\newcommand{\calC}{\mathcal{C}}
\newcommand{\calD}{\mathcal{D}}
\newcommand{\C}{\mathbb{C}}

\newcommand{\R}{\mathbb{R}}
\newcommand{\Z}{\mathbb{Z}}
\renewcommand{\P}{\mathbb{P}}

\title{Polynomial systems with few real zeroes}

\author{Benoit Bertrand}
\address{Departamento de \'Algebra\\
        Facultad de Ciencias Matem\'aticas\\
        Pza. de las Ciencias, 3\\
        Universidad Complutense de Madrid\\
        28040 Madrid\\ Spain}
\email{bertrand@mat.ucm.es}

\author{Frederic Bihan}
\address{Laboratoire de Math\'ematiques\\
         Universit\'e de Savoie\\
         73376 Le Bourget-du-Lac Cedex\\
         France}

\email{Frederic.Bihan@univ-savoie.fr}

\author{Frank Sottile}
\address{Department of Mathematics\\
         Texas A\&M University\\
         College Station\\
         TX \ 77843\\
         USA}
\email{sottile@math.tamu.edu}
\urladdr{http://www.math.tamu.edu/\~{}sottile}

\thanks{Part of work done at MSRI was supported by NSF grant DMS-9810361}
\thanks{Work of Sottile is supported by the Clay Mathematical Institute}
\thanks{Sottile and Bihan were supported in part by NSF CAREER grant DMS-0134860}
\thanks{Bertrand is supported by the European research network IHP-RAAG contract HPRN-CT-2001-00271}

\begin{document}

\begin{abstract}
  We study some systems of polynomials whose support lies in the
  convex hull of a circuit, giving a sharp upper bound for their numbers of
  real solutions.  This upper bound
  is non-trivial in that it is smaller than either the Kouchnirenko or
  the Khovanskii bounds for these systems.  When the support is
  exactly a circuit whose affine span is ${\Z}^n$, this bound is
  $2n+1$, while the Khovanskii bound is exponential in $n^2$.
  The bound $2n+1$ can be attained only
  for non-degenerate circuits.  Our methods 
  involve a mixture of combinatorics, geometry, and
  arithmetic.
\end{abstract}
\maketitle

%
\section*{Introduction}
The notion of degree for a multivariate Laurent polynomial $f$ is
captured by its support, which is the set of exponent vectors of monomials in
$f$.  For example, Kouchnirenko~\cite{Ko75} generalized the classical
B\'ezout Theorem, showing that the number of solutions in the complex
torus $(\mathbb{C}^*)^n$ to a generic system of $n$ polynomials in $n$
variables with common support $\calA$ is the volume $v(\calA)$ of the
convex hull of $\calA$, normalized so that the unit cube $[0,1]^n$ has
volume $n!$.  This Kouchnirenko number is also a (trivial) bound on
the number of real solutions to such a system of real polynomials.
This bound is reached, for example, when the common support $\calA$
admits a regular unimodular triangulation~\cite{St94b}.

Khovanskii gave a bound on the number of real solutions to a polynomial system
that depends only upon the cardinality $|\calA|$ of its support $\calA$:
 \[
   \mbox{Number of real solutions }\leq\  
    2^n 2^{\binom{|\calA|}{2}}\cdot (n+1)^{|\calA|}\,.
 \]
This enormous fewnomial bound is non-trivial (smaller than the Kouchnirenko
bound) only when the cardinality of $\calA$ is very small when compared to the
volume of its convex hull. 
It is widely believed that significantly smaller bounds should hold.
For example, Li, Rojas, and Wang~\cite{LRW03} showed that two trinomials in 2
variables have at most 20 common solutions, which is much less than the
Khovanskii bound of 995,328. 
While significantly lower bounds are expected, we know of no reasonable
conjectures about the nature of hypothetical lower bounds.

There are few other examples of non-trivial polynomial systems for
which it is known that not all solutions can be real via a bound
smaller than the Khovanskii bound.  We describe a class of supports
and prove an upper bound (which is often sharp) for the number of real
solutions to a polynomial system with those supports that is
non-trivial in that it is smaller than either the Kouchnirenko or
the Khovanskii bound.

A finite subset $\calA$ of ${\Z}^n$ that  affinely spans ${\Z}^n$ is 
{\it primitive}. 
A (possibly degenerate) {\it circuit} is a collection
$\calC:=\{0,w_0,w_1,\dotsc,w_n\}\subset\Z^n$ of 
$n{+}2$ integer vectors which spans $\R^n$. 
Here is the simplest version of our main results, which are proven in
Sections~\ref{Sec:upper} and~\ref{Sec:sharp}.\smallskip 

\noindent{\bf Theorem.}
{\it 
   A polynomial system with support a primitive circuit has at most
   $2n+1$ real solutions.
   There exist systems with support a primitive circuit having $2n+1$
   non-degenerate real solutions.
}\smallskip

This sharp bound for circuits suggests that
one may optimistically expect similar dramatic improvements in the doubly
exponential Khovanski bound for other sets $\calA$ of supports.

Adding the vectors $2w_0, 3w_0, \dotsc, kw_0$ to a circuit $\calC$ produces a
{\it near circuit} if $\R w_0\cap \calC = \{0,w_0\}$.  Let $\nu$ 
be the cardinality of
$\calD\cap\{w_1,\dotsc,w_n\}$, where $\calD\subset\calC$ is a minimal
affinely dependent subset.  Let $\ell$ be the largest integer so that
$w_0/\ell$ is integral.
We show that a polynomial system with support a primitive near circuit has at
most $k(2\nu{-}1)+2$ real solutions if $\ell$ is odd, or at most $2k\nu-1$
real solutions if $\ell$ is even, and these bounds are tight
among all such near circuits. These bounds coincide  when $k=1$, that
is, for primitive circuits.

For a given near circuit, we use its geometry and arithmetic to
give tighter upper bounds and construct systems
with many real zeroes.

An important step is to determine an eliminant for the system, which has the
form 
\[
     x^N \prod_{i=1}^p (g_i(x^\ell))^{\lambda_i}\ -\ 
         \prod_{i=p+1}^\nu (g_i(x^\ell))^{\lambda_i}\ ,
  \leqno{(*)}
\]
where the polynomials $g_i$ all have degree $k$, and $N$ and
$\lambda_i$ are the coefficients of the minimal linear dependence
relation among $\calD\cup\{w_0/\ell\}$.
 This reduces the problem to
studying the possible numbers of real zeroes of such a univariate
polynomial.  We adapt a method of Khovanskii to establish an upper
bound for the number of real zeroes of such a polynomial.  
Our upper bound uses a variant of the Viro construction, which also allows us to 
construct polynomials $(*)$ with many real zeroes.

In Section 1, we establish some basics on sparse polynomial systems,
and then devote Section 2 to an example of a family of near circuits
for which it is easy to establish sharp upper bounds, as simple linear
algebra suffices for the elimination, Descartes's rule of signs gives
the bounds, and the Viro construction establishes their sharpness.  In
Section 3, we compute an eliminant of the form $(*)$ for a system
supported on a near circuit.  Section 4 is devoted to proving an upper
bound for the number of real solutions to a polynomial of the form
$(*)$, while Section 5 constructs such polynomials with many real
zeroes, and some cases in which our bounds are sharp, including when $\calA$ is
a a circuit.

The authors wish to thank Stepan Orekov and Andrei Gabrielov for
useful discussions.

\section{Basics on sparse polynomial systems}\label{S:Basics}

We emphasize that we look for solutions to polynomial systems
which have only non-zero coordinates and thus lie in $(\C^*)^n$. 
Write $x^w$ for the monomial with exponent vector $w\in\Z^n$.
Let $\calA\subset\Z^n$ be a finite set 
which does not lie in an affine hyperplane.
A polynomial $f$ has {\it support} $\calA$ if 
the exponent vectors of its monomials lie in $\calA$.
A polynomial system with support $\calA$ is a system 
 \begin{equation}\label{E:system}
  f_1(x_1,\dotsc,x_n)\ =\  f_2(x_1,\dotsc,x_n)\ =\
   \dotsb \ =\ f_n(x_1,\dotsc,x_n)\ =\ 0\,,
 \end{equation}
where each polynomial $f_i$ has support $\calA$.
Multiplying a polynomial $f$ by a monomial $x^w$ does not change its set of
zeroes in $(\mathbb{C}^*)^n$, but does translate its support by the vector
$w$.
Thus it is no loss to assume that $0\in\calA$.
A system~\eqref{E:system} is {\it generic} if its number of solutions in
$(\C^*)^n$ equals $v(\calA)$, the volume of the convex hull of $\calA$
normalized so that the unit cube $[0,1]^n$ has
volume $n!$,
which is the Kouchnirenko bound~\cite{Ko75}.
This condition forces each solution to be simple.
We will always assume that our systems are generic in this sense.

Let $\Z\calA\subset \Z^n$ be the full rank sublattice generated by the vectors
in $\calA$.
(This is the affine span of $\calA$ as $0\in\calA$.)
If $\Z\calA=\Z^n$, then  $\calA$ is {\it primitive}.
The {\em index} of $\calA$ is the index of $\Z\calA$ in $\Z^n$.
The fundamental theorem of abelian groups implies that 
\[
   \frac{\Z^n}{\Z\calA}\ \simeq\
   \frac{\Z}{a_1\Z} \oplus \frac{\Z}{a_2\Z} \oplus\dotsb\oplus 
   \frac{\Z}{a_n\Z} \,,
\]
where $a_i$ divides $a_{i+1}$, for $i=1,\dotsc,n{-}1$.
These numbers $a_1,\dotsc,a_n$ are the {\it invariant factors of $\calA$}.
When $\calA$ is a simplex---so there are $n$
non-zero vectors in $\calA$---then the numbers $a_i$ are the invariant factors
of the matrix whose columns are these vectors.
The index of $\calA$ is the product of its invariant factors.

\subsection{Polynomial systems with support a simplex}
Let $e(\calA)$ be the number of even invariant factors of $\calA$.
The following result can be found in Section 3 of~\cite{St94b}.

\begin{prop}\label{P:Simplex}
 Suppose that $\calA$ is the set of vertices of a simplex.
 Then the number of real solutions to a generic system with support $\calA$ is
 \begin{itemize}
  \item[(i)]   $0$ or $2^{e(\calA)}$ if  $v(\calA)$ is even.
  \item[(ii)]  $1$ if $v(\calA)$ is odd.
 \end{itemize}
\end{prop}

\noindent{\it Proof.}
 Suppose that $0\in\calA$.
 Given a polynomial system~\eqref{E:system} with support $\calA$ whose
 coefficients are generic, we may perform Gaussian elimination on the 
 matrix of its coefficients and convert it into a system of the form
 \begin{equation}\label{E:Simplex}
   x^{w_i}\ =\ \beta_i\,,\quad\textrm{for}\quad i=1,\dotsc,n\,,
 \end{equation}
 where $\beta_i\neq 0$ and $w_1,\dotsc,w_n$ are the non-zero
 elements of $\calA$.  
 Solutions to this system have the
 form $\varphi^{-1}_\calA(\beta)$, where $\beta =
 (\beta_1,\dotsc,\beta_n) \in(\mathbb{C}^*)^n$
 and $\varphi_\calA$ is the homomorphism
\[
   (\mathbb{C}^*)^n\ \ni\ (x_1,\dotsc,x_n)\ \longmapsto\ 
   (x^{w_1},\dotsc,x^{w_n})\ \in\ (\mathbb{C}^*)^n\,.
\]
 Real solutions to~\eqref{E:Simplex} are $\psi^{-1}(\beta)$, where
 $\beta\in(\mathbb{R}^*)^n$ and
 $\psi\colon(\mathbb{R}^*)^n\to(\mathbb{R}^*)^n$ is the
 restriction of $\varphi_\calA$ to  $(\mathbb{R}^*)^n$.
 The kernel of $\psi$ consists of those points $x\in\{\pm1\}^n$ that satisfy
\[
  x^{w_i}\ =\ 1\quad \mbox{\rm for}\  i=1,\dotsc,n\,.
\]
 Let $A$ be a matrix whose columns are the non-zero elements of $\calA$.
 Identifying $\{\pm1\}$ with $\mathbb{Z}/2\mathbb{Z}$ identifies the
 kernel of $\psi$ with the kernel of the reduction of $A$ modulo 2,
 which has dimension $e(\calA)$.
 The result follows as $\psi$ is surjective if (and only if) 
 $e(A)=0$, which is equivalent to $v(\calA)$ being
 odd.
 \QED\vspace{-6pt}

\begin{rem}
 The upper bound of Proposition~\ref{P:Simplex} for a system with support
 $\calA$ is the volume of $\calA$ if and only if no invariant factor of
 $\calA$ exceeds 2.
\end{rem}

\subsection{Generic systems}
The system~\eqref{E:Simplex} is generic (has $v(\calA)$ simple roots in
$(\C^*)^n$) when the numbers $\beta_i$ are non-zero.
We give a proof of the following elementary result, as we will use the proof
later. 

\begin{prop}~\label{P:generic}
 Suppose that $\calA$ does not lie in an affine hyperplane.
 Then there is a non-empty Zariski open subset in the space of coefficients of 
 monomials appearing in a system with support $\calA$ such that the
 system has $v(\calA)$ simple solutions in $(\C^*)^n$.
\end{prop}

\noindent{\it Proof.}
 Suppose that $0\in\calA$.
 Then the affine span $\Z\calA$ is a full rank sublattice of $\Z^n$.
 Let $a:=|\calA|-1$ and consider the map $\varphi_\calA$ defined by
 \begin{equation}\label{E:mapPhi_A}
  \varphi_\calA\ \colon\  
   (\C^*)^n \ni (x_1,\dotsc,x_n)\ \longmapsto\
   [x^w\mid w\in\calA]\ \in\ \P^a\,.
 \end{equation}
 Its image is a subgroup of the dense torus in $\P^a$, and it is a
 homomorphism of algebraic groups.  Then a polynomial
 system~\eqref{E:system} with support $\calA$ is the pullback of $n$
 linear forms (given by the coefficients of the $f_i$) along the map
 $\varphi_\calA$.  These $n$ linear forms determine a linear section
 of the closure $X_\calA$ of the image of $\varphi_\calA$, a (not
 necessarily normal) projective toric variety~\cite[\S4,13]{St96}.
 Bertini's Theorem~\cite[p.~179]{Ha77} asserts that a general linear
 section of $X_\calA$ consists of $\deg(X_\calA)$ simple points, all
 lying in the image of $\varphi_\calA$.  Each of these pull back along
 $\varphi_\calA$ to $|\mbox{ker}(\varphi_\calA)|$ distinct solutions
 to the original system.

 The kernel of $\varphi_\calA$ is the abelian group dual to the factor
 group $\Z^n/\Z\calA$, which has order equal to the index of $\calA$.
 Furthermore, the toric variety $X_\calA$
 has degree equal to the volume of the convex hull of $\calA$, normalized so
 that a unit parallelepiped of $\Z\calA$ has volume $n!$.  
 The product of this volume with the index of $\calA$ is
 the usual volume of $\calA$.
 Thus a system with support $\calA$ and generic coefficients 
 has $v(\calA)$ simple solutions. \QED\vspace{-6pt}

\subsection{Congruences on the number of real solutions.}
 The proof of Proposition~\ref{P:generic} shows that the solutions to the 
 system~\eqref{E:system} are the fibres $\varphi_\calA^{-1}(\beta)$ for the
 points $\beta$ in a linear section of $X_\calA$.
 It follows that Proposition~\ref{P:Simplex} gives some restrictions on the
 possible numbers of real solutions. 
 The index of $\calA$ factors as  $2^{e(\calA)}\cdot N$.

\begin{cor}\label{C:triv_congr}
  The number of real solutions to~\eqref{E:system} is at most $v(\calA)/N$
  and is congruent to this number modulo $\max\{2,2^{e(\calA)}\}$.
\end{cor}

\noindent{\it Proof.}
Let $\varphi_\calA$ be the map~\eqref{E:mapPhi_A}.  The image of any
real solution to~\eqref{E:system} under $\varphi_\calA$ is a real
point, and there are $2^{e(\calA)}$ real solutions with the same image
under $\varphi_\calA$.  Thus the maximum number of real solutions
to~\eqref{E:system} is $2^{e(\calA)}\cdot\deg X_\calA$, which is
$v(\calA)/N$.  The congruence follows as the restriction of the map
$\varphi_\calA$ to the real subtorus is surjective on the real
subtorus of $X_\calA$ if and only if $e(\calA)=0$. \QED\vspace{-6pt}

\begin{rem}\label{R:Trivial}
 Any bound or
 construction for the number of real roots of polynomial systems associated to
 primitive vector configurations gives the same bounds and
 constructions for the number of real roots for configurations with odd
 index.
 Indeed, let $\calA\subset \Z^n$ be finite and $\calB\subset \Z^n$ be a basis for
 $\Z\calA$ so that linear combinations of vectors in $\calB$ identify $\Z\calA$ with
 $\Z^n$.
 This identification maps $\calA$ to a primitive vector configuration $\calA'$, which
 has the same geometry and arithmetic as $\calA$.
 The map $\varphi_\calA$~\eqref{E:mapPhi_A} factors
\[
  (\C^*)^n \ \xrightarrow{\, \varphi_{\calB\cup\{0\}}\,}\ 
  (\C^*)^n \ \xrightarrow{\, \varphi_{\calA'}\,}\ 
   \P^a\,. 
\]
 Since $\calA$ and $\calB\cup\{0\}$ have the same (odd) index,
 $\varphi_{\calB\cup\{0\}}$ is bijective on $(\R^*)^n$.
 Thus  $\varphi_{\calB\cup\{0\}}$ gives a bijection between real solutions to systems
 with support the primitive vector configuration $\calA'$ and real solutions to systems
 with support $\calA$.
\end{rem}

\section{A family of systems with a sharp bound}\label{S:exact}

We describe a family of supports $\Delta\subset\Z^n$ and prove a non-trivial 
sharp upper bound on the number of real solutions to polynomial systems with
support $\Delta$. 
That is, the points in $\Delta$ affinely span $\Z^n$,
but there are fewer than $v(\Delta)$ real solutions to polynomial systems with
support $\Delta$.
These sets $\Delta$ also have the property that they consist of all the integer
points in their convex hull.

Let $l > k>0$ and $n \ge 3$ be integers and 
$\epsilon = (\epsilon_1,\dotsc,\epsilon_{n-1})\in\{0,1\}^{n-1}$ be non-zero. 
Then
$\Delta_{k,l}^\epsilon \subset \R^n$ consists of the points
 \[
   (0,\dots,0),\ (1,0,\dots,0),\dots, (0,\dots,0,1,0),\ (0,\dots,0,k),\
   (\epsilon_1,\dots, \epsilon_{n-1},l)\,.
 \]
 together with the points
 along the last axis
\[
  (0,\dotsc,0,1),\ (0,\dotsc,0,2),\ \dotsc,\ (0,\dotsc,0,k{-}1)\,.
\]
Since these include the standard basis and the origin, $\Delta_{k,l}^\epsilon$
is primitive.

Set $|\epsilon|:=\sum_i\epsilon_i$.  Then the volume of
$\Delta_{k,l}^\epsilon$ is $l+k|\epsilon|$. Indeed,
the configuration $\Delta_{k,l}^\epsilon$ can be triangulated into
two simplices $\Delta_{k,l}^\epsilon \setminus \{(\epsilon_1,\dots, \epsilon_{n-1},l)\}$
and $\Delta_{k,l}^\epsilon \setminus \{0\}$ with volumes $k$ and $l-k+k|\epsilon|$,
respectively. One way to see this is to apply the affine transformation
\[ 
  (x_1,\dotsc,x_n)\ \longmapsto\ 
  (x_1,\dotsc,x_{n-1},x_n-k+k\sum_{i=1}^{n-1}x_i)\,.
\]

\begin{thm}
 The number, $r$, of real solutions to a generic system of $n$ real
 polynomials with support $\Delta_{k,l}^\epsilon$ lies in the interval
\[
    0 \ \leq\ r\ \leq \ k + k|\epsilon|+2\,,
\]
 and every number in this interval with the same parity as $l+k|\epsilon|$
 occurs. 
\end{thm}

This upper bound does not depend on $l$ and, since
$k < l$, it is smaller than or equal to the number $l+k|\epsilon|$
of complex solutions.
We use elimination to prove this result.

\begin{ex}\label{E:km-polytope}
Suppose that $n=k=3$, $l=5$, and $\epsilon = (1,1)$.  
\[
  \begin{array}{rclcl}
  \makebox[210pt][l]{Then the system}&&&&
   \multirow{9}*{
    \begin{picture}(80,105)(0,15)
      \ifEPSF 
        \put( 4,0){\epsfysize=4.5cm\epsffile{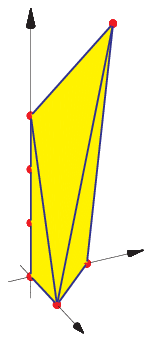}}
        \put(-1,24){$1$}   \put(16, 3){$x$}   \put(39,20){$y$}
        \put(-1,42){$z$}   \put(-1,63){$z^2$} \put(-1,85){$z^3$} 
        \put(50,120){$xyz^5$} 
      \else
        \put( 4,0){\fbox{\rule{0pt}{4.5cm}\rule{2cm}{0cm}}} 
      \fi
    \end{picture}}
  \\
    x + y + xyz^5 +1+z+z^2+z^3&=&0&\quad&\rule{0pt}{16pt}\\
    x + 2y + 3xyz^5 +5+7z+11z^2+13z^3&=&0&\quad&\\
   2x + 2y + xyz^5 +4+8z+16z^2+32z^3&=&0&\quad&\\
     \makebox[210pt][l]{is equivalent to }&&\rule{0pt}{16pt}\\
    x    -(5 +11z +23z^2 +41z^3)&=&0\rule{0pt}{16pt}\\
     y    +( 8 +18z +38z^2 +72z^3)&=&0\\
    xyz^5 -(2 +6z  +14 z^2 +30 z^3)&=&0\\
    \mbox{\ }
   \end{array}
\]
And thus its number of real roots equals the number of real roots of
\[
  z^5(5 +11z +23z^2 +41z^3)(8 +18z +38z^2 +72z^3) -
  (2 +6z  +14 z^2 +30 z^3)\,,
\]
%
%
which, as we invite the reader to check, is 3.
\end{ex}

\noindent{\it Proof.}
 A generic real polynomial system with support $\Delta_{k,l}^\epsilon$ has the
 form 
 \[
    \sum_{j=1}^{n-1} c_{ij}x_j \ + c_{in}x^\epsilon x_n^l\ \ 
    +\ f_i(x_n)\ =\ 0\ \quad 
    {\rm for\  }i=1,\dots,n\,,
 \]
 where each polynomial $f_i$ has degree $k$ and 
 $x^\epsilon$ is the monomial 
 $x_1^{\epsilon_1}\dotsb x_{n-1}^{\epsilon_{n-1}}$.
 

 Since all solutions to our system are simple, we may perturb the
 coefficient matrix $(c_{ij})_{i,j=1}^n$ if necessary and then use
 Gaussian elimination to obtain an equivalent system
 \begin{equation}\label{system}
  x_1-g_1(x_n)\ =\ \dotsb \ = \ x_{n-1}-g_{n-1}(x_n) \ =\ 
  x^\epsilon{x_n}^l -g_{n}(x_n)\ =\ 0\,,
 \end{equation}
 where each polynomial $g_i$ has degree $k$.
 Using the first $n-1$ polynomials to eliminate the
 variables $x_1, \dots , x_{n-1}$ gives the univariate polynomial
 \begin{equation}\label{E:eliminantex}
   {x_n}^l \cdot g_1(x_n)^{\epsilon_1} \dotsb 
    g_{n-1}(x_n)^{\epsilon_{n-1}}\ -\ g_n(x_n)\,,
 \end{equation}
 which has degree $l+k|\epsilon|=v(\Delta_{k,l}^\epsilon)$.  Any zero
 of this polynomial leads to a solution of the original
 system~\eqref{system} by back substitution.  This implies that the
 number of real roots of the polynomial~\eqref{E:eliminantex} is equal
 to the number of real solutions to our original
 system~\eqref{system}.

The eliminant~\eqref{E:eliminantex} has no terms of degree $m$ for $k<m<l$, and
so it has at most $k+k|\epsilon|+2$
non-zero real roots, by Descartes's rule of signs (\cite{De,BPR}, see
also Remark~\ref{L:Descartesrule}).
This proves the upper bound.  We complete the proof by constructing a
polynomial of the form~\eqref{E:eliminantex} having $r$ roots, for
every number $r$ in the interval between 0 and $k+k|\epsilon| + 2$
with the same parity as $l+k|\epsilon|$.

Choose real polynomials $f_1,\dotsc,f_n$ of degree $k$ having simple 
roots and non-zero constant terms.
We further assume that the roots of $f_1,\dotsc,f_{n-1}$ are distinct.
Put $f^\epsilon=f_1^{\epsilon_1}\dotsb f_{n-1}^{\epsilon_{n-1}}$ and
let $\alpha$ be the ratio of the leading term of $f_n$ to the constant term of 
$f^\epsilon$.
Choose any piecewise-linear convex function
$\nu:[0,k|\epsilon|+l]\rightarrow\R$ which is identically $0$ on
$[l,k|\epsilon|+l]$ and whose maximal domains of linearity 
are $[0,k]$, $[k,l]$, and $[l,k|\epsilon|+l]$.  

Let $f(z):= f_n \pm \alpha z^l f^\epsilon$ (the sign $\pm$ will
be determined later).  The {\it Viro polynomial} $f_t$ associated to $f$ and $\nu$
is obtained by multiplying the monomial $a_pz^p$ in $f(z)$ by $t^{\nu(p)}$
(\cite{Vi83,Vi84,Bi}).  By the definition of $\nu$,
$f_t=f_{n,t}+\alpha z^l f_1^{\epsilon_1}\dotsb
f_{n-1}^{\epsilon_{n-1}}$, where $f_{n,t}$ is the Viro polynomial
obtained from $f_n$ and the restriction of $\nu$ to $[0,k]$.  By
Viro's Theorem (\cite{Vi83,Vi84,Bi}, see also
Proposition~\ref{P:usefullemmaViro}), 
there exists a sufficiently small $t_0>0$, such that
if $t_0>t>0$, then the polynomial $f_t$ will have $r=r_1+r_2+r_3$
simple real roots, where $r_1$ is the number of real roots of
$f^\epsilon$, $r_2$ the number of real roots of $f_n$, and $r_3$ is
the number of (non-zero) real roots of the binomial obtained as the
truncation of $f$ to the interval $[k,l]$.  By our choice of
$\alpha$, this binomial is a constant multiple of $z^k \pm z^l$.  Thus
$r_3$ is $1$ if $l-k$ is odd, and either 0 or 2 (depending on the sign
$\pm$) if $l-k$ is even.  If $k$ is even, then every possible value of
$r$ between $0$ and $k+k|\epsilon| +2$ with the same parity as
$l+k|\epsilon|$ can be obtained in this way.  If $k$ is odd then this
construction gives all admissible values of $r$ in the interval
$[|\epsilon|+2,k+k|\epsilon| +2]$ but no values of  $r$ less than 
$|\epsilon|+2$.

Suppose now that $k$ is odd.  Take one of the Viro polynomials
$f_t$ with $r^\prime \geq |\epsilon|+2$ real roots.  For a generic
$t$, the crital values of $f_t$ are all different.  Choose one
such $t$ and suppose without loss of generality that the leading
coefficient of $f_t$ is positive.  For each $\lambda$, consider
the polynomial $h_\lambda= -\lambda - f_t$.  If $\lambda$ is
larger than every critical value of $f_t$, then $h_\lambda$ has
either 0 or 1 real roots, depending upon the parity of
$k|\epsilon|+l$.  Since the number of real roots of $h_\lambda$
changes by 2 when $\lambda$ passes through critical values of
$f_t$, and $h_0$ has at least $|\epsilon|+2$ real roots, every
possible number of real roots between $0$ and $|\epsilon|+2$ having
the same parity as $l+k|\epsilon|$ occurs for some $h_\lambda$.  \QED\vspace{-6pt}

\begin{rem}[Descartes's bound]\label{L:Descartesrule}
 Let $f(x)=\sum_{i=1}^d a_ix^{p_i}$ be a univariate polynomial with exponents 
 $p_1< \dotsb <p_d$ and $a_1,\dotsc,a_d$ are non-zero real numbers.
 Descartes's rule of signs asserts that the number of positive roots of $f$
 is no more than the number of $i \in\{1,\dotsc,d-1\}$ with 
 $a_i\cdot a_{i+1}<0$. 
 Applying this to $f(x)$ and $f(-x)$ shows that
 the number of non-zero real roots of $f$ is no more than $\sum_{i=1}^{d-1}
 \overline{(p_{i+1}-p_i)}$, 
 where $\overline{a}=1$ or $2$ according as $a$ is odd or even, respectively.
\end{rem}

\section{Elimination for near circuits}\label{S:elimination}

We first consider (possibly degenerate) circuits, which are collections of
$n{+}2$ integer vectors that affinely span $\R^n$.

\subsection{Arithmetic of circuits}\label{Sub:prelim}
Suppose that $\calC:=\{w_{-1},w_0,w_1,\dotsc, w_n\}\subset\Z^n$
affinely spans $\R^n$.  For each $i=-1,0,\dotsc,n$, let $\calA_i$ be
the (possibly degenerate) simplex with vertices
$\calC\setminus\{w_i\}$.  For any $j\in\{-1,0,\dotsc,n\}$, Cramer's
rule implies that
 \[
   \sum_{i=-1,\> i\neq j}^n (-1)^i \det(W_i) (w_i-w_j)\ \ =\ 0\,,
 \]
%
 where $W_i$ is the matrix whose columns are the vectors
 $w_{-1}-w_j,\dotsc,w_n-w_j$, with $w_i-w_j$ and $w_j-w_j$ omitted.
 Thus if $i\neq j$, $|\det W_i|$ is the volume $v(\calA_i)$ of
 $\calA_i$.

\begin{lemma}\label{L:primitive_circuit}
 Suppose that $\{0,v_0,\dotsc,v_n\}\subset\Z^n$ is primitive with 
 primitive relation
\[
  \sum_{i=0}^n\ \alpha_iv_i\ =\ 0\,.
\]
 If $\calA_q:=\{0,v_0,\dotsc,\widehat{v_q},\dotsc,v_n\}$, then 
 $\Z^n/\Z\calA_q \simeq \Z/\alpha_q\Z$.
\end{lemma}

\noindent{\it Proof.}
 We can assume that $q=0$.
 Since $\{0,v_0,\dotsc,v_n\}$ is primitive, $\Z^n=\Z\{v_0,v_1,\dotsc,v_n\}$.
 Then the image of $v_0$ generates the factor group 
 $\Z^n/\Z\{v_1,\dotsc,v_n\}$, and so this factor group is cyclic.
 Since the relation is primitive, $|\alpha_0|$ is the least positive
 multiple of $v_0$ lying in the lattice $\Z\{v_1,\dotsc,v_n\}$, which
 implies that  $\Z^n/\Z\{v_1,\dotsc,v_n\}\simeq \Z/\alpha_0\Z$. \QED

When $\calC$ is primitive, this proof shows that the invariant factors
of $\calA_q$ are $1,\dotsc,1, \alpha_q$, and $v(\calA_q)=\alpha_q$.
In general, the index $a$ of $\calC$ is the greatest common divisor of
the volumes $v(\calA_i)$, and the primitive affine relation on $\calC$
has the form
 \[
   \sum_{i=-1}^n \alpha_iw_i\ =\ 0\qquad\mbox{with}\qquad
    \sum_{i=-1}^n \alpha_i\ =\ 0\,,
 \]
 where $a|\alpha_i|=v(\calA_i)$.  Note that the configuration $\calC$
 admits two triangulations.  One is given by those $\calA_i$ with
 $\alpha_i>0$ and the other by those $\calA_i$ with
 $\alpha_i<0$.\smallskip

{}From now on, we make the following assumptions.  First, assume $w_{-1}=0$
and choose signs so that $\alpha_0\geq 0$.  Furthermore, assume the
vectors $w_1,\dotsc,w_n$ are ordered so that
$\alpha_1,\dotsc,\alpha_p>0$, $\alpha_{p+1},\dotsc,\alpha_\nu<0$, and
$\alpha_{\nu+1},\dotsc,\alpha_n=0$, for some integers $0\leq
p\leq\nu\leq n$.  If $a$ is the index of $\calC$ and we write
$\lambda_i=|\alpha_i|=v(\calA_i)/a$, then the primitive relation on
$\calC$ is
\[
   \sum_{i=0}^p\lambda_iw_i\ =\ \sum_{i=p+1}^\nu \lambda_iw_i\ ,
\]
and we have
 \begin{equation}\label{E:volumecircuit}
   v(\calC)\ =\ a\cdot \max\left\{
    \sum_{i=0}^p\lambda_i,\ \sum_{i=p+1}^\nu \lambda_i\right\}\,.
 \end{equation}

\subsection{Systems with support a near circuit}\label{nearcircuit}
Let $\{0,w_0,w_1,\dotsc,w_n\}\subset\Z^n$ span $\R^n$ and suppose that
$w_0=\ell e_n$, where $e_n$ is the $n$th standard basis vector and
$\ell$ is a positive integer.  Let $k>0$ and consider a generic
polynomial system with support
 \begin{equation}\label{E:Near_Circuit}
  \calC\ :=\ \{0,w_0, 2w_0,\dotsc, kw_0,\ w_1,\dotsc,w_n\}\,.
 \end{equation}
 We will call such a set of vectors a {\it near circuit} if no $w_i$
 for $i>0$ lies in $\R w_0$.  We will assume that $\calC$ is
 primitive, which implies that $\{0,e_n,w_1,\dotsc,w_n\}$ is
 primitive.  By Remark~\ref{R:Trivial}, this will enable us to deduce
 results for near circuits having odd index.

Write each vector $w_i=v_i+l_i\cdot e_n$, where $0\neq
v_i\in\Z^{n-1}$.  Then $\{0,v_1,\dotsc,v_n\}$ is primitive.  Let
$\lambda_1,\dotsc,\lambda_{\nu}$ be the positive integral coefficients
in the primitive relation on $\{0,v_1,\dotsc,v_n\}$,
 \[ 
   \sum_{i=1}^p \lambda_i v_i\ =\ 
   \sum_{i=p+1}^\nu \lambda_i v_i\ .
 \]
 Here, we could have $p=0$ or $p=\nu$, so that one of the two sums
 collapses to 0.

Assume that the vectors are ordered so that 
 \begin{equation}\label{E:primitiverelation}
   N e_n\ +\ \sum_{i=1}^p \lambda_i w_i \ -\ 
      \sum_{i=p+1}^\nu \lambda_i w_i\ =\ 0
 \end{equation}
is the primitive relation on $\{0,e_n, w_1,\dotsc,w_n\}$, where
 \[
   N\ :=\  \sum_{i=p+1}^\nu \lambda_il_i \ -\ 
          \sum_{i=1}^p \lambda_i l_i\ \geq\ 0\,.
 \]

 Consider a generic real polynomial system with support the near
 circuit $\calC$.  Perturbing the matrix of coefficients of the monomials $x^{w_i}$
 for $i=1,\dotsc,n$ and applying Gaussian elimination gives a system
 with the same number of real solutions, but of the form
 \begin{equation}\label{E:Csystem}
  x^{w_i}\ =\ g_i(x_n^\ell)\qquad \textrm{for}\ i=1,\dotsc,n\,,
 \end{equation}
where each $g_i$ is a generic polynomial of degree $k$.

Define a polynomial $f \in \R[x_n]$ by
 \begin{equation}\label{E:eliminant} 
   f(x_n)\ :=\  x_n^N \prod_{i=1}^p (g_i(x_n^\ell))^{\lambda_i}\ -\ 
         \prod_{i=p+1}^\nu (g_i(x_n^\ell))^{\lambda_i}\ .
  \end{equation}
Here, empty products are equal to 1.
By~\eqref{E:volumecircuit}, the degree of $f$ is equal to the volume of the
near circuit $\calC$, as $N$ is the volume of $\calA_0$ 
and $k\ell\lambda_i$ is the volume of the convex hull of 
$\{kw_0=k\ell e_n,0,w_1,\dotsc,\widehat{w_i},\dotsc,w_n\}$.
Lastly, the absolute value of the difference 
 \begin{equation}\label{E:Vol_last_simplex}
    \delta\ :=\ N +\ \sum_{i=1}^p k\ell\lambda_i  \ -\ 
      \sum_{i=p+1}^\nu k\ell\lambda_i 
 \end{equation}
in the degrees of the two terms of $f$ is the volume of the
convex hull of $\{k\ell e_n,w_1,\dotsc,w_n\}$, which may be zero.

\begin{thm}\label{T:strange_eliminant}
 Assume that $g_1,\dotsc,g_n$ are generic polynomials of degree $k$. 
 Then, the association of a solution $x$ of~$\eqref{E:Csystem}$ to its $n$th
 coordinate  $x_n$ gives a one-to-one correspondence between solutions
 of~$\eqref{E:Csystem}$ and roots of the univariate polynomial
 $f$~$\eqref{E:eliminant}$ which restricts to a bijection between real
 solutions to~$\eqref{E:Csystem}$ and real roots of $f$. 
\end{thm}

In particular, bounds on the number of real roots of polynomials of the
form~\eqref{E:eliminant} give bounds on the number of real solutions to
a generic polynomial system with support $\calC$.  
The polynomial $f$ is the eliminant of the system~\eqref{E:Csystem}, but our
proof is not as direct as the corresponding proof in Section~\ref{S:exact}.

\begin{lemma}\label{L:Two_Systems}
  For each $q=1,\dotsc,\nu$, the system 
 \[
  I_q\ \colon\ 
  \left\{\begin{array}{rcl} 
    x^{w_i}&=&g_i(x_n^\ell),\qquad\mbox{for }\ i=1,\dotsc,n, i\neq q\\
     f(x_n)&=&0
   \end{array}\right.
 \]
 is equivalent to the system
 \[
  J_q\ \colon\ 
  \left\{\begin{array}{rcl} 
    x^{w_i}&=&g_i(x_n^\ell),\qquad\mbox{for }\ i=1,\dotsc,n, i\neq q\\
    \big(x^{w_q}\big)^{\lambda_q}&=&\big(g_q(x_n^\ell)\big)^{\lambda_q}
   \end{array}\right.
 \]
\end{lemma}

\noindent{\it Proof.}
This follows from~\eqref{E:primitiverelation} and the
form~\eqref{E:eliminant} of $f$.  \QED\vspace{-6pt}

\begin{rem}\label{R:factor}
 Observe that the last polynomial in the system $J_q$ factors
 \[
   (x^{w_q})^{\lambda_q}\ -\ (g_q(x_n^\ell))^{\lambda_q}\ =\ 
   \prod_{\zeta\in Z_{\lambda_q}}\bigl(\zeta x^{w_q}\ -\ g_q(x_n^\ell)\bigr)\,,
 \]
 where $Z_{\lambda_q}$ is the set of roots of $z^{\lambda_q}-1$.
 For a solution $x$ to $J_q$, the number $\zeta$ is 
 $x^{-w_q}g_q(x_n^\ell)$.
 This factorization reveals that the system $J_q$ is a disjunction of
 $\lambda_q$ systems with support $\calC$.
 The subsystem with $\zeta=1$ is our original system~\eqref{E:Csystem}.
\end{rem}

\noindent{\it Proof of Theorem~$\ref{T:strange_eliminant}$.}
Since the integers $\lambda_1,\dotsc,\lambda_\nu$ are coprime, at
least one $\lambda_i$ is odd.  We restrict ourselves to the case where $\lambda_1$
is odd since the proof with any other $\lambda_i$ odd is similar.  Let
$q=1$ in Lemma~\ref{L:Two_Systems} and let $x_n$ be a root of $f$.  We
show that $x_n$ extends to a unique solution of~\eqref{E:Csystem}, and
that the solution is real if and only if $x_n$ is real.

 We first prolong $x_n$ to solutions to the system $I_1$ by solving the system
 for $y\in(\C^*)^{n-1}$
 \begin{equation}\label{E:Bsystem}
   y^{v_i}\ =\ x_n^{-l_i}g_i(x_n^\ell)\qquad \textrm{for}\ i=2,\dotsc,n\,.
 \end{equation}
 The numbers $\beta_i:=x_n^{-l_i}g_i(x_n^\ell)$ are well-defined and
 non-zero since the $g_i$ are generic polynomials of degree $k$.
 Hence ~\eqref{E:Bsystem} is a system associated to the simplex
 $\calB=\{0,v_2,\dotsc,v_n\}\subset\Z^{n-1}$ as studied in
 Proposition~\ref{P:Simplex}.  Its set of solutions is
 $\varphi_\calB^{-1}(\beta)$, where
 $\varphi_\calB\colon(\C^*)^{n-1}\to(\C^*)^{n-1}$ is the homomorphism
\[
  \varphi_\calB\ \colon\ t\ \longmapsto\ 
   (t^{v_2},\dotsc, t^{v_n})\ \in\ (\C^*)^{n-1}\,.
\]

 Let $y\in\varphi_\calB^{-1}(\beta)$.
 Then the fibre $\varphi_\calB^{-1}(\beta)$ is the set
 $\{ty\mid t\in \mbox{ker}(\varphi_\calB)\}$.
 We determine which of these, if any, is a solution to~\eqref{E:Csystem} by
 computing the number $\zeta$ of Remark~\ref{R:factor}.
 For $ty\in\varphi_\calB^{-1}(\beta)$, this number is
 $t^{-v_1} y^{-v_1}x_n^{-l_1}g_1(x_n^\ell)$.
 Consider the map 
 \[
     \psi \colon \mbox{ker}(\varphi_\calB) \ni t\  \longmapsto \ 
      t^{v_1} \in Z_{\lambda_1}.
\]
Since $\calD = \{0,v_1,\dotsc,v_n\}\subset\Z^{n-1}$ is primitive,
$\psi$ is well-defined by Lemma~\ref{L:primitive_circuit}.  Similarly,
$\varphi_\calD$ is injective (see Proposition~\ref{P:generic}), and
thus so is $\psi$.  Note that $\mbox{ker}(\varphi_\calB)$ is
isomorphic to $Z_{\lambda_1}$, thus $\psi$ is an isomorphism.  Thus
exactly one of these prolongations of $x_n$ given by
$\varphi_\calB^{-1}(\beta)$ is a solution to~\eqref{E:Csystem}.

 Suppose now that $x_n$ is real.
 Then $\beta$ is real.
 Since $\lambda_1$ is odd, there is a unique real solution
 $y$ to~\eqref{E:Bsystem}, by Proposition~\ref{P:Simplex}.
 But then $y^{-v_1}x_n^{-l_1}g_1(x_n^\ell)$ is real.
 As $\lambda_1$ is odd, there is eactly one real $\lambda_1$-th
root of unity, 
 namely 1, which proves that the real solution $(y,x_n)$ is a real solution
 to~\eqref{E:Csystem}. This completes the proof of the
 theorem. \QED\vspace{-6pt}

\begin{rem}\label{R:consprimitive}
  The primitivity of $\calC$ implies that $N$ and $\ell$ are coprime
  if $N \neq 0$, or $\ell=1$ if $N=0$.  Indeed, the affine span of
  $\calC$ is equal to that of $\{\ell e_n, w_1, \dotsc,w_n\}$ and the
  Cramer relation on this set is obtained by multiplying both sides
  of~\eqref{E:primitiverelation} by $\ell$. Hence, the index of
  $\calC$ is the greatest common divisor of $N,\ell \lambda_1, \dotsc,
  \ell \lambda_{\nu}$, and the result follows as the $\lambda_i$ are
  coprime. In particular, 
 either $\ell$ is odd, or $\ell$
  is even and $N$ and $\delta$ are odd.
\end{rem} 

\begin{ex}\label{E:construction_nearcircuit}
  We show that any positive integers $p$, $\ell$, $N$, with $N,\ell$
  coprime if $N \neq 0$ or $\ell=1$ if $N=0$, and any positive coprime
  integers $\lambda_1,\dotsc,\lambda_\nu$ with $\nu\leq n$ correspond
  to a primitive near circuit, when one $\lambda_i=1$.  Thus any
  polynomial of the form~\eqref{E:eliminant} is the eliminant of a
  system with support a primitive near circuit, when one exponent
  $\lambda_i=1$.

Assume without loss of generality that $\lambda_\nu=1$.
Let $e_1,\dotsc,e_n$ be the standard basis in ${\R}^n$. 
Let $v_i:=e_i$ for $i=1,\dotsc,\nu-1$.
Set 
\[
  v_\nu\ :=\ \sum_{i=1}^p \lambda_i v_i\ -\  \sum_{i=p+1}^{\nu-1} \lambda_i v_i\,.
\]
Since $\lambda_1,\dotsc,\lambda_\nu$ are coprime, there exist integers
$l_1,\dotsc,l_\nu$ such that $N=\sum_{i=p+1}^\nu\lambda_il_i-
\sum_{i=1}^p\lambda_il_i$.  If we set $w_i:=v_i+l_ie_n$ for
$i=1,\dotsc, \nu$
 and $w_i:=e_{i-1}$ for $i=\nu+1,\dotsc,n$, then we
obtain the relation~\eqref{E:primitiverelation} among
$e_n,w_1,\dotsc,w_n$.  It is then easy to see that the near circuit
$\{0,\ell e_n,\dotsc,k\ell e_n,w_1,\dotsc,w_n\}$ is primitive (for any
integer $k$) if we assume that $N$ and $\ell$ are coprime, or $\ell=1$
if $N=0$.
\end{ex}

\section{Upper bounds for near circuits}\label{Sec:upper}

We first give a version of Viro's construction for univariate polynomials that
takes multiplicities into account.
We then use this to establish upper bounds for the number of real roots of
polynomials of the form~\eqref{E:eliminant} by studying the total variation in
the number of real roots of a pertubation of the eliminant. 

\subsection{Viro univariate polynomials}

Consider a univariate Viro polynomial
 \[
    f_t(y)\ =\ \sum_{p=p_0}^d \phi_p(t) \, y^p \,,
 \]
where $t$ is a positive real number, and each coefficient
$\phi_p(t)$ is a finite sum $\sum_{q \in I_p} c_{p,q}t^{q}$ with
$c_{p,q}\in\R$ and $q$ a rational number. 
Write $f$ for the function of $y$ and $t$ defined by $f_t$.

Let $P$ be the convex hull of the points $(p, q)$ for $p_0 \leq p \leq d$
and $q \in I_p$. Assume that $P$ has dimension $2$.
Its lower hull $L$ is the union of the edges $L_1,\dotsc,L_l$ of $P$ whose
inner normals have positive second coordinate. 
Let $I_i$ be the image of $L_i$ under the projection to the first axis.
Then the intervals $I_1,\dotsc,I_l$ subdivide the Newton segment $[p_0,d]$ of
$f_t$. 

Let $f^{(i)}$ be the facial subpolynomial of $f$ for the face $L_i$.
That is, $f^{(i)}$ is the sum of terms  $c_{p,q}y^p$ such that $(p,q) \in L_i$.
Suppose that $L_i$ is the graph of $y \mapsto a_iy+b_i$ over $I_i$,
Expanding $f_t(y t^{-a_i}) / t^{b_i}$ in powers of $t$ gives
 \begin{equation}\label{E:changecoord}
  \frac{f_t(y t^{-a_i})}{t^{b_i}}\ =\
    f^{(i)}(y)+ t^{A_i}d^{(i)}(y) + h^{(i)}(y,t) \;  , \quad i=1,\dotsc,l,
 \end{equation}
 where $t^{A_i}d^{(i)}(y)$ collects the terms with smallest positive
 power of $t$ and $h^{(i)}(y,t)$ collects the remaining terms (whose
 powers of $t$ exceed $A_i$).  Then $f^{(i)}(y)$ has Newton segment
 $I_i$ and its number of non-zero roots counted with multiplicities is
 $|I_i|$, the length of the interval $I_i$.

\begin{prop}\label{P:usefullemmaViro}
 Assume that for any non-zero root $\rho$ of $f^{(i)}$, $i=1,\dotsc,l$,
 either $\rho$ is a simple root of $f^{(i)}$, or else $d^{(i)}(\rho)\neq 0$.
 Then there exists $t_0>0$ such for $0<t<t_0$, the univariate polynomial
 $f_t(y)$ has only simple non-zero roots, with 
\[ 
  r\ =\ \sum \; c(\rho)
\]
 non-zero real roots, where the sum is over $i=1,\dotsc,l$ and then all
 non-zero real roots $\rho$ of $f^{(i)}$ where 
 \[ 
   c(\rho)\ =\ \left\{
   \begin{array}{ll}
    1 & \mbox{if the multiplicity $m$ of $\rho$ is odd},\\
    0 & \mbox{if $m$ is even and $f^{(i)}(y)/d^{(i)}(\rho)>0$, for $y$ near $\rho$,}\\
    2 & \mbox{if $m$ is even and $f^{(i)}(y)/d^{(i)}(\rho)<0$, for $y$ near $\rho$.}
   \end{array}\right.
 \]
\end{prop}

In particular, if the non-zero roots of $f^{(1)}, \dotsc, f^{(l)}$ are
simple, then the number of non-zero real roots of $f_t$ for $t>0$
small enough equals the total number of non-zero real roots of
$f^{(1)}, \dotsc, f^{(l)}$.  This is the usual version of Viro's
theorem for univariate polynomials.  \medskip

\noindent{\it Proof.}
 For each root $\rho \neq 0$ of $f^{(i)}(y)$ of multiplicity $m$, there will be
 $m$ roots near $\rho$ to
\[
   f^{(i)}(y)+ t^{A_i}d^{(i)}(y) + h^{(i)}(y,t)\,,
\]
 for $t>0$ sufficiently small.
 This gives $|I_i|$ roots to $f_t(y t^{-a_i})/t^{b_i}$, and thus 
 all solutions to $f_t$ in $\C^*$, at least when $t>0$ is
 sufficiently small. 
 Indeed, let $K\subset\C^*$ be a compact set containing the non-zero roots of
 the facial polynomials $f^{(1)}(y),\dotsc,f^{(l)}(y)$.
 Then, for $t>0$ sufficiently small, $K$ contains the $|I_i|$ roots to 
 $f_t(y t^{-a_i})/t^{b_i}$ that we just constructed.
 The compact sets $t^{-a_1}K,\dotsc,t^{-a_l}K$ are pairwise disjoint for
 $t>0$ sufficiently small,
 and this gives $|I_1|+\dotsb+|I_l|=d-p_0$ non-zero simple roots of $f_t$ for
 $t>0$ small enough. 
 But this accounts for all the non-zero simple roots of $f_t$.
 
 We now determine how many 
 roots of $f_t(y t^{-a_i})/t^{b_i}$ are real.  
 Roots close to $\rho$ are real only if $\rho$ is real, 
 and then the number of such real roots is determined by the first
 two terms $f^{(i)}(y)+ t^{A_i}d^{(i)}(y)$ in $t$, as
 $d^{(i)}(\rho)\neq 0$.  But this polynomial has $c(\rho)$ real roots
 near $\rho$.   \QED\vspace{-6pt}

\subsection{Upper bounds}\label{S:Upper_Bounds}

We give upper bounds on the number of real roots of a generic polynomial system
with support a primitive near circuit
 \[
  \calC\ =\ \{0,\ell e_n, 2\ell e_n,\dotsc, k\ell e_n,\ w_1,\dotsc,w_n\}\,.
 \]
As explained in Section~\ref{S:elimination}, it suffices to bound the roots of 
a polynomial $f$ of the form~\eqref{E:eliminant}.
Consider the polynomial $f_t(y)$ depending on a real parameter
$t\neq 0$ defined by  
 \[
   f_t(y)\ :=\ t \cdot y^N \prod_{i=1}^p (g_i(y^\ell))^{\lambda_i}\ -\ 
         \prod_{i=p+1}^\nu (g_i(y^\ell))^{\lambda_i}\ ,
 \]
where $g_1,\dotsc,g_\nu$ are generic polynomials of degree $k$.
We will study how the number of real roots of $f_t$ can vary as $t$ runs from 
$\infty$ to $1$.
Note that $f_1$ is our original eliminant $f$.

\begin{rem}\label{R:elimsyst}
If $p \neq 0$ then $f_t(x_n)$ is the eliminant of the system
 \[
    \left\{ 
    \begin{array}{lll}
         x^{w_1}  & = & t^{1/\lambda_1} \cdot g_1(x_n^\ell)\,, \\
         x^{w_i} & = & g_i(x_n^\ell) \quad i=2,\ldots,n\,.
    \end{array}\right.
 \]
If $p=0$, then $f_t$ is $t^2$ times the eliminant of the system
 \[
   \left\{ 
    \begin{array}{lll}
         x^{w_1} & = & t^{-1/\lambda_1} \cdot g_1(x_n^\ell)\,,\\
         x^{w_i} & = & g_i(x_n^\ell) \quad i=2,\ldots,\ldots,n\,.
     \end{array}\right.
 \]
\end{rem}

For any integer $a$, define $\overline{a} \in \{0,1,2\}$ by 
\[
  \overline{a}\ :=\  \; \left\{\begin{array}{ll}
     2 & \mbox{if $a$ is positive and even}\\
     1 & \mbox{if $a$ is positive and odd}\\
     0 & \mbox{otherwise}
      \end{array}\right.
\]
A root $\rho$ of a univariate polynomial is {\em singular} if it has 
multiplicity greater than 1. 

Let $\chi(Y)$ be the boolean truth value of $Y$, so that 
$\chi(0>1)=0$, but $\chi(0<1)=1$.

\begin{prop}\label{P:totalmultroots}
  The total sum of the multiplicities of the non-zero singular real roots
  of $f_t$ for $t \in \R^*$ is no more than $2k\overline{\ell}\nu -
  2\bar{\ell}\left(\chi(\delta=0)+\chi(N=0)\right)$.
  
  Moreover, if $\ell$ is even, so that $N$ is odd due to the
  primitivity of $\calC$, then the total sum of the multiplicities of
  the non-zero singular real roots of $f_t$ for $t>0$ is equal to the
  corresponding sum for $t<0$. Hence, both numbers are no more
  $2k\nu$.
\end{prop}

\noindent{\it Proof.}
Write $f_t=tF-G$, where $F$ and $G$ are the two terms of $f$.
Let $\rho$ be a non-zero root of $f_t$ for some $t \neq 0$.
Then $F(\rho)G(\rho)\neq 0$, as the roots of
$g_1,\dotsc,g_\nu$ are distinct. 
Note that $t=G(\rho)/F(\rho)$.
Then $\rho$ is a singular root of $f_t$ if and only if
\[
   (F'G-FG')(\rho)\ =\ 0\,.
\]
If $N\neq 0$, then the polynomial $F'G-FG'$ factors as
 \begin{equation}\label{factor}
   (F'G-FG')(y)\ =\ 
     \left(y^{N-1} \prod_{i=1}^\nu (g_i(y^\ell))^{\lambda_i-1}\right) \cdot H(y)
 \end{equation}
where $H$ is the polynomial defined by
 \[
  H(y)\ =\ \prod_{i=1}^\nu g_i(y^\ell) \cdot
   \left(N+ \ell y^\ell\cdot D(y^\ell) \right)\,,
 \]
with 
\[
   D(z)\ =\ \sum_{i=1}^p \lambda_i \cdot \frac{g'_i(z)}{g_i(z)}\ -\ 
           \sum_{i=p+1}^\nu \lambda_i \cdot \frac{g'_i(z)}{g_i(z)}\,.
\]
(If $N=0$, then $y^{N-1}$ is replaced by $y^{\ell-1}$, and the last factor in
$H$ is simply $\ell D(y^\ell)$.) 

Thus $H(y)=h(y^\ell)$, where $h$ is a polynomial of degree $k\nu -
\left(\chi(\delta=0)+\chi(N=0)\right)$ with a non-zero constant term,
as the $g_i$ are generic.  If $\rho$ is a non-zero singular real root
of $f_t$ for $t \neq 0$ then $h(\rho^\ell)=0$.  Thus the total number
of non-zero singular real roots of $f_t$ for $t \in\R^*$ is at most
$k\overline{\ell}\nu-\bar{\ell}\left(\chi(\delta=0)+\chi(N=0)\right)$.

A root $\rho$ of $f_t$ has multiplicity $m \geq 2$ for some $t \neq 0$
if and only if $tF^{(i)}(\rho)-G^{(i)}(\rho)=0$ for $i=0,\dotsc,m-1$
and $tF^{(m)}(\rho)-G^{(m)}(\rho) \neq 0$.
This is equivalent to the system
 \begin{equation}\label{multipleforf}
  \begin{array}{rcll}
   F^{(i)}(\rho) \cdot  G^{(j)}(\rho) -F^{(j)}(\rho)  \cdot G^{(i)}(\rho)
   & =  &0& \mbox{\quad if} \quad 0 \leq i,j \leq m-1, \\
   F^{(i)}(\rho) \cdot G^{(j)}(\rho)-F^{(j)}(\rho) \cdot G^{(i)}(\rho)
   & \neq &0& \mbox{\quad if $i=m$ and $0 \leq j \leq m-1$}.\rule{0pt}{13pt}
  \end{array}
 \end{equation}
Solving~\eqref{factor} for $H$ and using the expression of the $k$th derivative
of $F'G-FG'$ as a sum of polynomials of the form
$F^{(i)}G^{(j)}-F^{(j)}G^{(i)}$, 
we can use~\eqref{multipleforf} to deduce that $\rho$ is a root of multiplicity
$m\geq 2$ of $f_t$ for some $t \neq 0$ if and only if $\rho$ is a root of
multiplicity $m{-}1$ of $H$. 
Thus the total sum of the multiplicities of the non-zero singular real
roots of $f_t$ for $t\in\R^*$ is 
%
 \begin{equation}\label{boundmultipleH}
   \sum_{\mbox{$\rho$ a real root of $H$}}(m_{\rho}(H)+1)\,,
 \end{equation}
 where $m_{\rho}(H)$ is the multiplicity of the root $\rho$ of $H$.
 We see that~\eqref{boundmultipleH} is bounded by
 $2k\overline{\ell}\nu-2
 \overline{\ell}\left(\chi(\delta=0)+\chi(N=0)\right)$ with equality
 when all roots of $h$ are real and simple (and positive if $\ell$ is even) so
 that the singular real roots of $f_t$ for $t\in\R^*$ are real double roots.

Finally, the statement concerning the case $\ell$ even and $N$ odd is
obvious after noting that in this case the (non-zero) real roots of
$h(y^{\ell})$ come in pairs $(\rho,-\rho)$, the function $G/F$ is an
odd function, and $N$ and $\delta$ are both odd integer numbers. 
\QED\vspace{-6pt}

\begin{rem}\label{R:delta}
  If $\{k\ell e_n,w_1,\dotsc,w_n\}$ are affinely independent and
  $N\neq 0$, then the difference $\delta$~\eqref{E:Vol_last_simplex}
  of the degrees of the terms of $f$ is non-zero, and the polynomials
  $f_t$ have the same Newton segment for $t\neq 0$.  Thus the number
  of real roots of $f_t$ can change only if $t$ passes through $0$, or
  through a value $c \neq 0$ such that $f_c$ has a singular root.
 
 If the difference $\delta=0$, then there is one number $t_\infty$ for
 which the degree of $f_{t_\infty}$ drops.  If necessary, we may
 perturb coefficients of one $g_i$ so that the number of real roots of
 $f$ does not change, and the degree of $f_{t_\infty}$ drops by one.
 This will result in no net change in the number of real roots of
 $f_t$ as $t$ passes through $t_\infty$, for the root which
 `disappears' in $f_{t_\infty}$ is a real root at infinity.
 Similarly, if $N=0$, then we may assume that there is one number $t_0$ for which the  
 constant term of $f_{t_0}$ vanishes.  Perturbing again if necessary
 results in no net change in the number of non-zero real roots of
 $f_t$ as $t$ passes through $t_0$.
 
 Thus, the number of values $c$ where the number of real roots of
 $f_t$ changes is finite by Proposition~\ref{P:totalmultroots}, and
 hence it makes sense to define the numbers
 \[
   r_{-\infty} \; , \quad r_{0-} \; , \quad r_{0+} \; , \quad r_{+\infty}
 \]
 as the numbers of real roots of $f_t$ as $t$ tends to $-\infty$, $0$
 by negative values, $0$ by positive values and $+\infty$,
 respectively.
\end{rem}

Recall that $\chi(Y)$ denotes the boolean truth value of $Y$.

\begin{prop}\label{P:extremalrealroots}
We have
 \begin{eqnarray*}
   \frac{r_{0+}+r_{0-}}{2} & \leq & k\overline{\ell}(\nu-p)+ \chi(\delta >0)\\
   \frac{r_{+\infty}+r_{-\infty}}{2} & \leq & 
        k\overline{\ell}p+\chi(N>0)+\chi(\delta < 0)\\
   \frac{|r_{0+}-r_{0-}|}{2} & \leq & 
    k\overline{\ell} \left(\sum_{i=p+1}^\nu \overline{\lambda_i}\right) 
      -k\overline{\ell}(\nu-p)+
      \chi(\delta>0 \mbox{ is even})\\
    \frac{|r_{+\infty}-r_{-\infty}|}{2} & \leq  & 
     k\overline{\ell} \left(\sum_{i=1}^p \overline{\lambda_i}\right) 
      -k\overline{\ell}p+
      \chi(N>0\mbox{ is even}) + \chi(\delta<0\mbox{ is even})
 \end{eqnarray*}
Furthermore, if $l$ is even and $N$ is odd, we have
 \begin{eqnarray*}
\frac{r_{0+}+r_{+\infty}}{2} & \leq & k \nu +1
 \end{eqnarray*}
\end{prop}

\noindent{\it Proof.}
 As in the proof of Proposition~\ref{P:totalmultroots}, write
 $f_t=tF-G$.
 We apply Proposition~\ref{P:usefullemmaViro} (and its proof) to $f_t$ and
 $f_{-t}=-tF-G$ to estimate $r_{0+}$ and $r_{0-}$, respectively.

 Let $P$ be the common Newton polygon of $f_t(y)$ and $f_{-t}(y)$, as
 polynomials in $t$ and $y$.
 Projecting the lower faces of $P$ onto the first coordinate axis gives a 
 single interval $I_1$ if $\delta\leq 0$, or the union of two intervals $I_1$
 and $I_2$ if $\delta>0$.
 Here, $I_1=[0,\deg(G)]$, the Newton segment of $G$ and 
 $I_2=[\deg(G),\deg(F)]$, which has length $\delta$. 
 For both $f_t(y)$ and $f_{-t}(y)$,
 the polynomial $f^{(1)}$ corresponding to $I_1$ is just $G$.
 If $\delta>0$, then the polynomial $f^{(2)}$ corresponding to $I_2$
 is the binomial $\pm M_F-M_G$, which is the difference of the highest degree
 terms of $\pm F$ and $G$.

 Both binomials $\pm M_F-M_G$ have only simple non-zero roots.
 The polynomial $f^{(1)}=G$ has singular roots if any of 
 $\lambda_{p+1},\dotsc,\lambda_\nu$ are not equal to $1$. 
 Since $F$ and $G$ have no common root, the assumptions of
 Proposition~\ref{P:usefullemmaViro} are fulfilled for both $f_t$ and $f_{-t}$.
 The numbers $r_{0+}$ and $r_{0-}$, hence $r_{0+}+r_{0-}$ and
 $|r_{0+}-r_{0-}|$, are sums of contributions of the non-zero real roots of
 $G$ and of the non-zero real roots of $\pm M_F-M_G$.
 
 Consider contributions from roots $\rho$ of $G$, which satisfy
 $\rho^\ell=\zeta$, where $\zeta$ is a root of some $g_i$ for
 $p<i\leq\nu$.  This has multiplicity $\lambda_i$.  If $\lambda_i$ is
 odd, then $\rho$ contributes 1 to both $r_{0\pm}$ and hence 2 to
 $r_{0+}+r_{0-}$ and 0 to $|r_{0+}-r_{0-}|$.  If $\lambda_i$ is even,
 then $\rho$ contributes 2 or 0 to $r_{0+}$, depending upon the sign
 of $G/F$ near $\rho$.  Replacing $t$ by $-t$, shows that it
 contributes 2 or 0 to $r_{0-}$, depending upon the sign of $-G/F$
 near $\rho$.  Thus $\rho$ contributes 2 both to $r_{0+}+r_{0-}$ and
 to $|r_{0+}-r_{0-}|$.
 
 Suppose now that $\delta>0$.  Then each non-zero real root of the
 binomial $\pm M_F-M_G$ is simple and thus contributes 1 to
 $r_{0\pm}$.  Both binomials have only one non-zero real root if
 $\delta$ is odd.  If $\delta$ is even, then $M_F-M_G$ (resp.
 $-M_F-M_G$) has $0$ or $2$ (resp. $2$ or $0$) real roots according as
 the product of the coefficients of $M_F$ and $M_G$ is positive or
 negative.  It follows that the roots of these binomials contribute
 $2$ to $r_{0+}+r_{0-}$, and contribute $0$ or $2$ to
 $|r_{0+}-r_{0-}|$ according as $\delta$ is odd or even, respectively.
 Summing up all contributions gives the desired upper bounds for
 $r_{0+}+r_{0-}$ and $|r_{0+}-r_{0-}|$.
 
 The upper bounds for $r_{+\infty}-r_{-\infty}$ and
 $|r_{+\infty}-r_{-\infty}|$ are obtained in exactly the same way if
 we use the polynomial $g_t(y)=F(y)-tG(y)$ instead of $f_t(y)$.  The
 Newton polygon $Q$ of $g_t$ is the reflection of $P$ in the
 horizontal line of height $1/2$, so that the lower faces of $Q$ are
 the upper faces of $P$.  They project to 1, 2, or 3 intervals on the
 first coordinate axis, $I_1=[0,N]$ (if $N\neq 0$), $I_2=[N,\deg(F)]$,
 and $I_3=[\deg(F),\deg(G)]$, if $\delta<0$.  The polynomial $f^{(2)}$
 corresponding to $I_2$ is just $F$, and the other polynomials are
 binomials.

 Finally, assume that $\ell$ is even, $N$ is odd, and let us prove the
 last inequality.  Using the facts that the non-zero real roots of $G$
 come in pairs $(\rho,-\rho)$, the function $G/F$ is an odd function,
 and that $N$ and $\delta$ are both odd integers, we obtain $r_{0+}
 \leq 2k (\nu - p )+ \chi(\delta >0)$.  Similarly, using the
 polynomial $g_t(y)=F(y)-tG(y)$ instead of $f_t(y)$, we obtain that
 $r_{+\infty} \leq 2kp+1+\chi(\delta<0)$. Suming up these two
 inequalities gives the result.  \QED\vspace{-6pt}

\begin{thm}\label{T:upperboundunivariate}
 The number $r$ of real solutions to a generic system with support the near
 circuit $\calC$ satisfies the following inequalities
 \begin{equation}\label{E:firstboundunivariate}
 \begin{array}{lll}
  r & \leq & 2k\overline{\ell}p
          + k\overline{\ell} \left(\sum_{i=p+1}^\nu \overline{\lambda_i}\right) 
        +\chi(N>0) + 1 - \chi(\delta>0\mbox{ is odd})\\
  &  & \\
 & & -\chi(\delta=0)-\overline{\ell} \left(\chi(\delta=0)+\chi(N=0)\right) \,,
\end{array}
 \end{equation}
and
 \begin{equation}\label{E:secundboundunivariate}
 \begin{array}{lll}
r & \leq &  2k\overline{\ell}(\nu-p)
       + k \overline{\ell}\left(\sum_{i=1}^p \overline{\lambda_i}\right)
    + \chi(N>0\mbox{ is even}) + 1-\chi(\delta<0\mbox{ is odd}) \\
  &  & \\
 & &
-\chi(\delta=0)-\overline{\ell} \left(\chi(\delta=0)+\chi(N=0)\right)\,,
\end{array}
 \end{equation}
 where $k,\ell,N,\nu$, and $\lambda_i$ are defined in
 Section~\ref{nearcircuit}.

Moreover, if $\ell$ is even and $N$ is odd, then we have
 \begin{equation}\label{E:thirdboundunivariate}
r \leq 2 k \nu + 1
\end{equation}
\end{thm}

\begin{rem}\label{R:absolutebound}
 Since $\lambda_1,\dotsc,\lambda_\nu$ are relatively prime, we see that the
 absolute upper bound for such a near circuit with $\ell$ odd is
\[
   k(2\nu-1)+2\,,
\]
and this can be obtained if exactly one $\lambda_i$ is odd, $N>0$ and
even, and $\delta$ is even and non-zero.  This upper bound is
maximized when $\{0,e_n,w_1,\dotsc,w_n\}$ forms a non degenerate
circuit, that is, if no proper subset is affinely dependent.  If
$\ell$ is even, the corresponding absolute upper bound is $2 k \nu +
1$, and is also maximized when $\{0,e_n,w_1,\dotsc,w_n\}$ forms a non
degenerate circuit.
\end{rem}

\noindent{\it Proof of Theorem~$\ref{T:upperboundunivariate}$.}
Let $f$ be the univariate eliminant of a generic polynomial system
with support $\calC$, which has the form~\eqref{E:eliminant}.  For an
interval $I\subset\R$, let $\Delta_{I}$ be the (positive) variation in
the number of non-zero real roots of $f_t$ for $t \in I$.  As noted in
Remark~\ref{R:delta}, the number of non-zero real roots of $f_t$ can
change only if $t$ passes through $0$, or if $t$ passes through a
value $c \neq 0$ such that $f_c$ has a real singular root.  Passing
through the value $t=0$, the variation of the number of real roots of
$f_t$ is at most $|r_{0+}-r_{0-}|$.

 Recall that $f=f_t$ for $t=1$ and that $f$ has no singular roots.
 Considering the path from $t$ close to $-\infty$ to $t=1$, we obtain
 \[
   r\ \leq\ r_{-\infty}+\Delta_{(-\infty,0)}+|r_{0+}-r_{0-}|+\Delta_{(0,1)}\,.
 \]
 Considering the path from $t$ close to $+\infty$ to $t=1$, we obtain
 \[
   r\ \leq\ r_{+\infty} + \Delta_{(1,+\infty)}\,.
 \]
 Combining these two inequalities yields
 \begin{equation}\label{E:bound1}
   r\ \leq\ \frac{r_{+\infty}+r_{-\infty} +|r_{0+}-r_{0-}|}{2}
     +\frac{\Delta_{\R^*}}{2}\,.
 \end{equation}
 The number $\Delta_{\R^*}$ is at most the total multiplicity of the
 singular real roots of $f_t$ for $t \neq 0$.
  By
 Proposition~\ref{P:totalmultroots}, this is at most
 $2k\overline{\ell}\nu-2
 \overline{\ell}\left(\chi(\delta=0)+\chi(N=0)\right)$.  The
 inequality~\eqref{E:firstboundunivariate} follows then
 from~\eqref{E:bound1} and Proposition~\ref{P:extremalrealroots}.

 Using the polynomial $g_t=F-tG$ in place of $f_t=tF-G$,
 leads to
 \begin{equation}\label{E:bound2}
  r\ \leq\ \frac{r_{0+}+r_{0-} +|r_{+\infty}-r_{-\infty}|}{2}
        +\frac{\Delta_{\R^*}}{2}.
 \end{equation}
 The inequality~\eqref{E:secundboundunivariate} is then obtained using
 Proposition~\ref{P:extremalrealroots}.

 Finally, considering the paths from $t=0$ to $t=1$, and from $t$
 close to $+\infty$ to $t=1$ gives
\begin{equation}\label{E:bound3}
  r\ \leq\ \frac{r_{0+}+r_{+\infty}}{2} + \frac{\Delta_{(0,+\infty)}}{2}.
 \end{equation}
 For $\ell$ even and $N$ odd, the inequality
 ~\eqref{E:thirdboundunivariate} comes then from the corresponding
 statements in Proposition~\ref{P:totalmultroots} and
 Proposition~\ref{P:extremalrealroots}. \QED\vspace{-6pt}

\section{Constructions and sharp upper bounds}\label{Sec:sharp}

We now construct polynomials $f$ having the form~\eqref{E:eliminant} with many
real roots.  In some cases, this 
achieves the upper bound
of Theorem~\ref{T:upperboundunivariate} for the maximal number of real
solutions to generic systems with support a given near circuit.

\begin{thm}\label{T:lowerboundpol}
 Let $\calC$ be a primitive near circuit with $k$, $\ell$, $N$, and $\lambda_i$
 as in Section~$\ref{nearcircuit}$.
 Suppose that $d_1, \ldots, d_\nu$ are nonnegative
          integers with $d_i\leq k$ such that
    \begin{equation}\label{E:condlotroots}
         \ell \sum_{i=1}^\nu d_i \lambda_i\ <\ N+k\ell \sum_{i=1}^p \lambda_i\,.
    \end{equation}
If $\ell$ is odd,
     then there is a generic polynomial system with support $\calC$ having 
  \[
      \sum_{i=1}^\nu  d_i \,\overline{\lambda_i} \; + \; 
       \overline{d}
  \]
     real solutions, where $d$ is the (positive) difference of the two sides
     of~\eqref{E:condlotroots}.
If $\ell$ is even, hence $N$ is odd,
 then there is a generic polynomial system with support $\calC$ having 
  \[
 2 \sum_{i=1}^{\nu}d_i + 1
  \]
     real solutions.
\end{thm}

We use Proposition~\ref{P:usefullemmaViro} to determine the number of real
roots for small $t>0$ of a polynomial of the form
 \[
   f_t(x)\ =\ 
     t^a \prod_{\zeta \in R_1}(t^{-b}x^{\ell}-\zeta)^{m(\zeta)} \ - \ 
    x^{\mu}\prod_{\zeta \in R_2}(\zeta-x^{\ell})^{m(\zeta)}\,,
 \]
where $\mu,a,b$ are positive integers, $R_1$ and $R_2$ are disjoint sets of
positive real numbers, and $m(\zeta)$ is a positive integer for 
$\zeta\in R_1\cup R_2$. 
Set $\mu_i:=\ell \cdot \sum_{\zeta\in R_i} m(\zeta)$.

\begin{lemma}\label{L:mainlemlowerbound}
 Suppose that $\mu_1<a \ell /b<\mu$, and we have further that
\begin{enumerate} 
 \item 
      if $\zeta,\zeta' \in R_1$ with $m(\zeta)$ even and $m(\zeta')$ odd, 
      then $\zeta'< \zeta$, and 
 \item 
      if $\zeta,\zeta' \in R_2$ with $m(\zeta)$ even and $m(\zeta')$ odd, 
      then $\zeta'> \zeta$.
\end{enumerate}

Let $e$ (respectively $o$) be the number of $\zeta\in R_1\cup R_2$
such that $m(\zeta)$ is even (respectively odd).  If $\ell$ is odd,
then, for $t>0$ sufficiently small, $f_t$ has exactly
$2e+o+\overline{\mu - \mu_1}$ simple non-zero real roots.  If $\ell$
is even and $\mu$ is odd, then, for $t>0$ sufficiently small, $f_t$
has exactly $2e+2o+1$ simple non-zero real roots.
\end{lemma}

\noindent{\it Proof.}
 We use the notation of Proposition~\ref{P:usefullemmaViro}.
 The inequalities $\mu_1 < a\ell /b < \mu$ imply that the lower hull of $P$
 consists of three segments whose projection onto the first coordinate axis are
 the intervals $I_1=[0,\mu_1]$, $I_2=[\mu_1,\mu]$, and $I_3=[\mu,\mu+\mu_2]$. 
 The corresponding facial subpolynomials are
 \begin{eqnarray*}
  f^{(1)}(x)&=& \prod_{\zeta \in R_1}(x^{\ell}-\zeta)^{m(\zeta)}\\
  f^{(2)}(x)&=& x^{\mu_1}\;-\;
      \big(\prod_{\zeta \in R_2}\zeta^{m(\zeta)}\big)\cdot x^\mu\\
  f^{(3)}(x)&=& -x^\mu \prod_{\zeta \in R_2}(\zeta-x^{\ell})^{m(\zeta)}\,.
 \end{eqnarray*}
 Note that $\prod_{\zeta \in R_2}\zeta^{m(\zeta)}>0$ since $ R_2$ consists of
 positive real numbers. 
 Thus the non-zero roots of the binomial $f^{(2)}$ are simple, with 
 $\overline{\mu-\mu_1}$ of them real.
 Proposition~\ref{P:usefullemmaViro} applies as we can see from the
 expansions~\eqref{E:changecoord} for $f^{(1)}$ and $f^{(3)}$.
 \begin{eqnarray*}
   f_t(xt^{\frac{b}{\ell}})/t^a &=&
    f^{(1)}(x)\;-\;t^{\frac{b\mu}{\ell} -a}x^{\mu}
          \prod_{\zeta \in R_2}(\zeta-x^{\ell}t^b)^{m(\zeta)}\\
    &=&f^{(1)}(x)\;-\; t^{\frac{b\mu}{\ell}-a}x^\mu(\prod_{\zeta\in R_2}\zeta^{m(\zeta)})
      \; +\; h^{(1)}(x,t)\,.\\
   f_t(x) &=& 
    f^{(3)}(x)\; + \; t^a\prod_{\zeta \in R_1}(t^{-b}x^{\ell}-\zeta)^{m(\zeta)}\\
     &=& f^{(3)}(x)\; + \; t^{a-\frac{b\mu_1}{\ell}} x^{\mu_1} \;+\; h^{(3)}(x,t)\,.
 \end{eqnarray*}
 Assume that $\ell$ is odd. Then the map $x \mapsto x^{\ell}$ is a
 bijection from the real roots of $f^{(1)}$ (resp., $f^{(3)}$) to
 $R_1$ (resp., $R_2$).  Conditions $(1)$ and $(2)$ imply that the
 contribution of any real root $\rho$ such that $\rho^{\ell}=\zeta \in
 R_i$ and $m(\zeta)$ is even is equal to $2$.
 
 Assume that $\ell$ is even and $\mu$ is odd. Then the real roots of
 $f^{(1)}$ (resp., $f^{(3)}$) come in pairs $(\rho,-\rho)$ with
 $\rho^{\ell}=\zeta \in R_1$ (resp., $R_2$).  If $m(\zeta)$ is odd,
 then the contributions of $\rho$ and $-\rho$ are both equal to $1$.
 If $m(\zeta)$ is even, then one 
 contribution is $2$, while the other is $0$, 
 as $\mu$ is odd. 
 Finally, 
 note that $\overline{\mu-\mu_1}=1$ if $\ell$ is even and $\mu$ is
 odd.  \QED

\noindent{\it Proof of Theorem~$\ref{T:lowerboundpol}$.}
 Set
\[
  \mu\ :=\ N+\ell\sum_{i=1}^p (k-d_i)\lambda_i\,.
\]
 The inequality~\eqref{E:condlotroots} can be rewritten as
\[ 
  \ell\sum_{i=p+1}^\nu  d_i\lambda_i\ <\  \mu\,.
\]

Then, by Lemma~\ref{L:mainlemlowerbound} there exist polynomials
$h_1,\dotsc,h_n$ with distinct roots such that $h_i$ has degree $d_i$
with $d_i$ real roots and the polynomial
 \[
  g(x)\ :=\  x^\mu \prod_{i=1}^p (h_i(x^\ell))^{\lambda_i}\ -\ 
         \prod_{i=p+1}^\nu (h_i(x^\ell))^{\lambda_i}
 \]
has either 
 \[
  \sum_{i=1}^\nu d_i\, \overline{\lambda_i}\; + \; \overline{d}
 \]
or else
\[ 2 \sum_{i=1}^{\nu}d_i + 1\]
simple real roots according as $\ell$ is odd, or $\ell$ is even and
$N$ is odd, respectively.  The polynomial $g(x)$ can be rewritten as
 \[
   g(x)\ =\ 
     x^N \prod_{i=1}^p \bigl(x^{\ell(k-d_i)}h_i(x^\ell)\bigr)^{\lambda_i}\ -\ 
         \prod_{i=p+1}^\nu \bigl(h_i(x)^\ell\bigr)^{\lambda_i}
 \]

 If $d_i=k$, set $g_i(x):=h_i(x)$. 
 Otherwise, set 
\begin{eqnarray*}
  g_i(x)\ :=\ \epsilon(1+x+\cdots+x^{k-d_i-1})+x^{k-d_i}h_i(x)
    \qquad 1 \leq i \leq p\\
  g_i(x)\ :=\ h_i(x)+\epsilon(x^{d_i+1}+\cdots+x^k)
    \qquad p+1 \leq i \leq \nu\,.
\end{eqnarray*}
 For sufficiently small $\epsilon >0$, the polynomial
 \[
   f(x)\ =\ x^N \prod_{i=1}^p \bigl(g_i(x^{\ell})\bigr)^{\lambda_i}-  
    \prod_{i=p+1}^\nu \bigl(g_i(x^{\ell})\bigr)^{\lambda_i}
 \]
 has simple roots and at least the same number of real roots as $g$. \QED\vspace{-6pt} 

\begin{thm}\label{T:sharpbounds1}\rule{0pt}{20pt}
 Assume that $N > k\ell\sum_{i=p+1}^\nu \lambda_i$ and let $m$ be the maximal
 number of real solutions to a generic system with support the near circuit
 $\calC$. 
\smallskip

If $\ell$ is even, then  $m=2k\nu+1$.\smallskip

Suppose now that $\ell$ is odd. 
\begin{enumerate}
\item If $\lambda_1,\dotsc,\lambda_p$ are even, then
  \[ 
    m\ =\ 2kp 
      +k\sum_{i=p+1}^\nu\overline{\lambda_i}+\overline{\delta}\,.
  \]

\item If exactly one number among $\lambda_1,\dotsc,\lambda_p$ is odd,
  $k=\ell=1$ and $\delta$ is odd, then
 \[
   m\ =\ 2p+1+\sum_{i=p+1}^\nu \overline{\lambda_i}\,.
 \]

\item If $\lambda_{p+1},\dotsc,\lambda_n$ are even, then
 \[
   m\ =\ 2k(n-p)+k\sum_{i=1}^p
   \overline{\lambda_i} + \overline{N}\,. 
 \]

\item If exactly one number among $\lambda_{p+1},\dotsc,\lambda_n$ is odd,
  $k=\ell=1$ and $N$ is odd, then
 \[
   m\ =\ 2(n-p)+1+\sum_{i=1}^p \overline{\lambda_i}\,.
 \]
\end{enumerate}

\end{thm}

\noindent{\it Proof.}
We apply Theorem~\ref{T:lowerboundpol} with each $d_i=k$.  The case of 
$\ell$ even is a direct consequence of
Theorem~\ref{T:upperboundunivariate} and
Theorem~\ref{T:lowerboundpol}.  Suppose now that $\ell$ is odd and let
$B_1$ and $B_2$ be the upper bounds for the number of real solutions
to a generic system with support $\calC$ which are given in
Theorem~\ref{T:upperboundunivariate} by
formulas~\eqref{E:firstboundunivariate}
and~\eqref{E:secundboundunivariate}, respectively.  Set
\[
  B_-\ :=\ k\sum_{i=1}^\nu \overline{\lambda_i} + \overline{d}\,,
 \]
 where $d:=N-k \ell \sum_{i=p+1}^\nu \lambda_i >0$. Note that
 $\delta,N \geq d$, so $\delta, N>0$.  By
 Theorem~\ref{T:lowerboundpol}, the number $B_-$ is a lower bound on
 the maximal number of real solutions of a generic system with support
 $\calC$.

We check that $B_1\geq B_-$ and analyze the conditions under which $B_1=B_-$.
As $\delta>0$ and $N >0$, we have 
\[
   B_1-B_-\ =\ k\sum_{i=1}^p (2-\overline{\lambda_i})\, 
          + \, \overline{\delta}-\overline{d}\,.
\]
We have
\[
  \overline{\delta}-\overline{d}\ =\ \ 
\left\{
\begin{array}{rl}
0 & \mbox{ if $k\ell\sum_{i=1}^p \lambda_i$ is even} \\
1 & \mbox{ if $k\ell\sum_{i=1}^p \lambda_i$ is odd and $d$ is odd}\rule{0pt}{13pt}\\
-1 & \mbox{if $k\ell\sum_{i=1}^p \lambda_i$ is odd and $d$ is even}\rule{0pt}{13pt}
 \end{array}
\right.
\]
If $\overline{\delta}-\overline{d}=0$, then $B_1 \geq B_-$ with
equality only if $\lambda_1,\dotsc,\lambda_p$ are even.  This proves
Part $(1)$.

If $\overline{\delta}-\overline{d}=1$, then $B_1 > B_-$.
Assume now that $\overline{\delta}-\overline{d}=-1$. 
Then $k\ell \sum_{i=1}^p \lambda_i$ is odd, $d$ is even
and $B_1-B_-=k\ell\sum_{i=1}^p (2-\overline{\lambda_i}) -1$. 
Since $k\ell\sum_{i=1}^p \lambda_i$ is odd, $k\ell$ is odd and at least one
number among $\lambda_1,\dotsc,\lambda_p$ is odd. 
Thus $B_1-B_- \geq k\ell-1$ with equality only if exactly one number among
$\lambda_1,\dotsc,\lambda_p$ is odd. 
Part $(2)$ now follows.

Parts $(3)$ and $(4)$ are similar.  \QED\vspace{-6pt} 

\begin{thm}\label{T:coeff1or2}
  Assume that $\lambda_i \in \{1,2\}$ for $i=1,\dotsc,\nu$, $\ell$ is
  odd, and let $m$ be the maximal number of real solutions to a
  generic system with support the near circuit $\calC$.

\begin{enumerate}
 \item If $N >  k \ell \sum_{i=p+1}^\nu \lambda_i$,
 then 
 \[
   m\ =\ k\sum_{i=1}^\nu \lambda_i + 
    \overline{{\textstyle N-k\ell \sum_{i=p+1}^\nu \lambda_i}}\,.
 \]

\item Suppose that $\ell=1$.
     If $N< k \sum_{i=p+1}^\nu \lambda_i$, then 
 \[
   m\ =\ v(\calC)\ = \  \max \left\{ k \sum_{i=p+1}^\nu \lambda_i \ ,
                      N+k\sum_{i=1}^p \lambda_i\right\}\ .
\]
\end{enumerate}
\end{thm}

\noindent{\it Proof.}
 For Part (1), the number $m$ equals the upper bound given by Descartes's
 rule of signs when applied to a polynomial of the form~\eqref{E:eliminant}. 
 Theorem~\ref{T:lowerboundpol} with each $d_i=k$ implies the existence of a
 polynomial with this form with $m$ real roots.

 For Part (2), we have that  $N< k \sum_{i=p+1}^\nu \lambda_i$.
 If we also have $N+k\sum_{i=1}^p \lambda_i > k\sum_{i=p+1}^\nu \lambda_i$,
 then there exist nonnegative integers $d_1,\dotsc,d_\nu \leq k$ such that
\[
   d\ =\ N+k\sum_{i=1}^p \lambda_i- \left(\sum_{i=1}^\nu \lambda_i d_i\right)
    \  \in\ \{1,2\}\,,
\]
as $\lambda_i\in\{1,2\}$.  By Theorem~\ref{T:lowerboundpol}, there
exists a polynomial of the form~\eqref{E:eliminant} having
$\sum_{i=1}^{\nu} \lambda_i d_i + d=N+k\sum_{i=1}^p\lambda_i$ non-zero
real roots.
 
 Finally, suppose that $ N+k\sum_{i=1}^p\lambda_i\le k \sum_{i=p+1}^\nu
 \lambda_i$. 
 Consider a polynomial
\[ 
  f_t(x)\ =\ t \cdot x^N \prod_{i=1}^p (g_i(x))^{\lambda_i}-
  \prod_{i=p+1}^\nu (g_i(x))^{\lambda_i}\,,
\]
where $g_1,\dotsc,g_\nu$ are polynomials of degree $k$ with non-zero
constant terms.  The lower part of the Newton polygon of $f_t(x)$
consists of a single segment projecting onto $[0, k \sum_{i=p+1}^\nu
\lambda_i]$.  Hence, if Proposition~\ref{P:usefullemmaViro} applies,
the number of real roots of $f_t$ for $t>0$ small enough is the sum of
contributions of the non-zero real roots of $g_{p+1},\dotsc,g_\nu$.
Choosing polynomials $g_1,\dotsc,g_\nu$ with distinct roots satisfying
conditions
\begin{enumerate}
\item If $i,j>p$, then $g_i$ has $k$ positive roots, and if
  $\lambda_i$ is odd and $\lambda_j$ is even, then the leading
  coefficient of $g_i$ is positive and every root of $g_i$ is less
  than every root of $g_j$.
\item 
   The polynomials $g_1,\dotsc,g_p$ are positive at each root of 
   of $g_{p+1},\dotsc, g_\nu$.
\end{enumerate}
 By Proposition~\ref{P:usefullemmaViro}, $f_t$ has
 $k\sum_{i=p+1}^\nu\lambda_i$ non-zero real roots for $t>0$ small enough.
  \QED\vspace{-6pt}  
 
\begin{rem}
 The example of Section~2 is a special case of Part(1) of
 Theorem~\ref{T:coeff1or2}.
 Indeed, in Section~2, we have
\[
   w_0=e_n,\ w_i=e_i\ i=1,\dotsc,n{-}1,\ \mbox{and}\ 
   w_n=\sum \epsilon_i e_i\ +\ le_n\,,
\]
 so that 
\[
  N\ =\ l, \ \nu\ =\ |\epsilon|+1\ =\ \sum_i\lambda_i\ ,
  \ \mbox{and}\ p\ =\ |\epsilon|\,.
\]
 Then the maximum number of Part(1) of
 Theorem~\ref{T:coeff1or2} is
\[
  m\ =\ k(|\epsilon|+1) + \overline{l-k}\,,
\]
 which is what we found in Section 1.
\end{rem}

\begin{thm}\label{T:sharpabsolutebound}
  The number of real roots of a generic system with support a
  primitive near circuit $\calC=\{0,\ell e_n,\dotsc,k\ell
  e_n,w_1,\dotsc,w_n\}$ in ${\R}^n$ is at most $(2\nu-1)k+2$ if $\ell$
  is odd, or $2k\nu+1$ if $\ell$ is even.  Moreover, these bounds are
  sharp.
\end{thm}

\noindent{\it Proof.}
As the numbers $\lambda_1,\dotsc,\lambda_\nu$ are coprime, at least
one is odd.  Since $\overline{\lambda_j} \leq 2$, the upper bound for
odd $\ell$ follows from Theorem~\ref{T:upperboundunivariate} (see
Remark~\ref{R:absolutebound}).  The sharpness of this bound follows
from Theorem~\ref{T:coeff1or2} for a primitive near circuit with all
$\lambda_1,\dotsc,\lambda_\nu$ but one equal to $2$ and one which is
equal to $1$, and where and $N-k{\ell}\sum_{i=p+1}^\nu \lambda_i$ is
positive and even (Example~\ref{E:construction_nearcircuit} shows that
such a near circuit exists).
 
 The bound $2k\nu+1$ for $\ell$ even comes from
 Theorem~\ref{T:upperboundunivariate}, its sharpness follows from
 Theorem~\ref{T:sharpbounds1}. \QED

A near circuit with $k=1$ is just a circuit. The following result is a
particular case of the previous one.

\begin{thm}
The number of real roots of a generic system with support a primitive
circuit $\calC=\{0,\ell e_n,,w_1,\dotsc,w_n\}$ in
 ${\R}^n$ is at most $2\nu+1$, and this bound this sharp.
 
 The absolute upper bound for the number of real roots of a generic
 system with support a primitive circuit in ${\R}^n$ is $2n+1$, and
 this bound is sharp.  Moreover, this bound can be attained only for
 non-degenerate circuits.
\end{thm}

\noindent{\it Proof.}
We only need to prove the last sentence as the others are 
corollaries of Theorem~\ref{T:sharpabsolutebound}.  For this, we note
that the bound $2\nu+1$ is obtained when $\delta$ and $N$ are
non-zero.  Hence the absolute bound $2n+1$ is obtained when $\nu=n$,
$\delta$ and $N$ are non-zero, which is exactly the case of a
non-degenerate circuit.   \QED\vspace{-6pt}

\providecommand{\bysame}{\leavevmode\hbox to3em{\hrulefill}\thinspace}
\providecommand{\MR}{\relax\ifhmode\unskip\space\fi MR }
\providecommand{\MRhref}[2]{%
  \href{http://www.ams.org/mathscinet-getitem?mr=#1}{#2}
}
\providecommand{\href}[2]{#2}

\end{document}

\bibliographystyle{amsplain}
\bibliography{bibl}
 
\end{document}